\documentclass{cocv}

\usepackage[applemac]{inputenc}   
\usepackage[T1]{fontenc}
\usepackage{lmodern}
\usepackage{graphicx}
\usepackage[english]{babel}
\usepackage{amsmath, amsfonts, amssymb, amsthm}
\usepackage{hyperref} 
\hypersetup{colorlinks = true, allcolors=blue}
\usepackage{array}
\usepackage{enumitem}
\makeatletter 
\renewcommand\theequation{\thesection.\arabic{equation}} 
\@addtoreset{equation}{section}
\makeatother

\newcommand{\vertiii}[1]{{\left\vert\kern-0.25ex\left\vert\kern-0.25ex\left\vert #1 
    \right\vert\kern-0.25ex\right\vert\kern-0.25ex\right\vert}}
\newcommand{\changed}[1]{{\color{black}#1}}

\begin{document}

\title{State-constrained controllability of linear reaction-diffusion systems }
\author{Pierre Lissy}\address{CEREMADE,  Universit\'e Paris-Dauphine \& CNRS UMR 7534, Universit\'e PSL, 75016 Paris, France;\\ \email{lissy@ceremade.dauphine.fr}}
\author{Cl\'{e}ment Moreau}\sameaddress{1}
\date{\today}
\runningtitle{Necessary conditions for local controllability...}

\begin{abstract}
We study the controllability of a coupled system of linear parabolic equations, with nonnegativity constraint on the state. We establish two results of controllability to trajectories in large time: one for diagonal diffusion matrices  with an ``approximate'' nonnegativity constraint, and a another stronger one, with ``exact'' nonnegativity constraint, when all the diffusion coefficients are equal and the eigenvalues of the coupling matrix have nonnegative real part. The proofs are based on a ``staircase'' method. Finally, we show that state-constrained controllability admits a positive minimal time, even with weaker unilateral constraint on the state.
\end{abstract}

\begin{resume}
On s'int\'eresse \`a la contr\^olabilit\'e d'un syst\`eme coupl\'e d'\'equations paraboliques lin\'eaires avec une contrainte de positivit\'e sur l'\'etat. On \'enonce deux r\'esultats de contr\^olabilit\'e aux trajectoires en temps grand : un pour des matrices de diffusion diagonales avec contrainte de positivit\'e ``approch\'ee'', et un autre, plus fort, avec une contrainte de positivit\'e ``exacte'', dans le cas o\`u les coefficients de diffusion sont identiques et o\`u les valeurs propres de la matrice de couplage sont de partie r\'eelle positive. Les preuves s'appuient sur une m\'ethode ``en escalier''. Enfin, on montre que le temps minimal de contr\^olabilit\'e avec contrainte sur l'\'etat est strictement positif, y compris sous une contrainte unilat\'erale moins restrictive sur l'\'etat. 
\end{resume}




\maketitle

\textit{This article is dedicated to Enrique Zuazua for the occasion of his $60$th birthday, with admiration for his outstanding achievements and all the new paths he explores in control theory, that are sources of inspiration notably for the younger generation.}

\section{Introduction}


\changed{In the following, $\mathbb{N}$ and $\mathbb{N}^*$ denote the sets of respectively nonnegative and positive integers.} Let $d$ in $\mathbb{N}^*$. Let $\Omega$ be a bounded open connected set of $\mathbb{R}^d$ with $\mathrm{C}^\infty$ boundary $\partial \Omega$, and $\omega$ an open subset of $\Omega$. Let $T>0$, $\Omega_T = (0,T) \times \Omega$ and $\omega_T = (0,T) \times \omega$. Let $n, m$ in $\mathbb{N}^*$, with $n\geqslant 2$. Let $\overrightarrow{n}$ be the outward normal on $\partial\Omega$. We consider the following parabolic linear system of $n$ coupled scalar equations with homogeneous Neumann boundary conditions and internal control:

\begin{equation}
\label{eq:paralin}
\left \{
\begin{array}{l l l}
\partial_t Y - D \Delta Y &= A Y + B U \mathbf{1}_{\omega} &\text{in } \Omega_T, \\
\partial_{\overrightarrow{n}} Y&= 0 &\text{on } (0,T) \times \partial \Omega, \\
Y(0,\cdot)&=Y^0 (\cdot) &\text{in } \Omega.
\end{array}
\right .
\end{equation}

In \eqref{eq:paralin}, $A = (a_{ij})_{1\leqslant i,j \leqslant n}$ and $D = (d_{ij})_{1\leqslant i,j \leqslant n}$ are square matrices in $\mathcal{M}_n(\mathbb{R})$ and $B$ is a matrix in $\mathcal{M}_{n,m}(\mathbb{R})$. Assume that $D$ satisfies the ellipticity condition given by
\begin{equation}
\exists \alpha >0, \forall \xi \in \mathbb{R}^n, \langle D \xi, \xi \rangle \geqslant \alpha | \xi |^2.
\label{eq:ellipticity}
\end{equation}

The spaces in which  the initial condition $Y^0$ and the control $U$ lie will be made more precise later on. Notice that $m$ represents the number of controls. Notably, we may have $m<n$, which means that we can have an underactuated system.

The controllability  to trajectories of system \eqref{eq:paralin} in arbitrary time has been established under a Kalman-type condition in \cite{ammar-khodja_kalman_2008}. 
The question we address in this paper is the following : is it possible to ensure that the state remains nonnegative while controlling \eqref{eq:paralin} from a nonnegative initial state towards a nonnegative trajectory? This question is relevant because reaction-diffusion systems like \eqref{eq:paralin} frequently model phenomena in which the state is nonnegative (\textit{e.g.} concentrations of chemicals). In these cases, a controlled trajectory that does not remain nonnegative would have no interest for applications. 

State-constrained controllability is a challenging subject that has gained popularity in the last few years, notably at the instigation of J\'erome Loh\'eac, Emmanuel Tr\'elat and Enrique Zuazua in the seminal paper \cite{loheac_minimal_2017}, in which some controllability results with positivity constraints on the state or the control for the linear heat equation are proved, under a minimal time condition which turns out to be necessary. This question yielded to an increasing number of articles in different frameworks, many of them being coauthored by Enrique Zuazua: for ODE systems \cite{loheac_minimal_2018,loheac2020nonnegative},  semilinear and quasilinear heat equations \cite{pighin_controllability_2018,nunez2019controllability}, monostable and bistable reaction-diffusion equations \cite{mazari2020constrained,ruiz2020control}, the fractional one-dimensional heat equation \cite{antil2019controllability,biccari2019controllability}, wave equations \cite{pighin2019controllability}, and age-structured systems \cite{maity2019controllability}. The spirit of most of these results can be summarized this way: when the considered system is controllable in the classical sense, and sometimes under assumptions on the initial and target states or on the system properties, controllability with a constraint on the state is possible with a positive minimal time.

Our goal in this paper is to state similar results for coupled parabolic systems of the form \eqref{eq:paralin} with internal control. This framework raises several difficulties. Indeed, for boundary control or equations satisfying a maximum principle, there is an equivalence between nonnegativity of the state and nonnegativity of the control, which is useful as control-constrained problems are better understood in general. This equivalence does not hold anymore for System \eqref{eq:paralin}: the state might remain positive even if the control is negative. Moreover, due to the coupling terms, the asymptotic behaviour of the trajectories is difficult to know precisely. This means that we cannot rely on dissipativity or stabilization to a steady state to obtain controllability, with the notable exception of the case where the diffusion matrix $D$ is equal to the identity matrix. 

One of the ideas we present in the following to bypass this difficulty is an original version of the ``staircase'' method, where we drive the system along a path of non-constant trajectories, while this approach is usually employed to follow a path of constant steady states \cite{loheac_minimal_2017,pighin_controllability_2018}. Adapting it to non-constant trajectories requires supplementary arguments to make sure that these trajectories do not go too far away from each other.

In the particular case where $D=I_n$ and the eigenvalues of $A$ have nonnegative real part, we establish controllability to trajectories in large time with nonnegative state (Theorem \ref{thrm:main-id}). When $D$ is only assumed to be diagonal, we show that \eqref{eq:paralin} is also controllable in large time, but with state remaining ``approximately'' nonnegative, \textit{i.e.} greater than $-\varepsilon$ for any fixed $\varepsilon>0$ (Theorem \ref{thrm:main}). Additionally, we show that there exists a positive minimal controllability time as soon as the initial state and the target trajectory are different, even if we allow the state to be greater than a negative constant instead of being nonnegative (Theorem \ref{thrm:minimal}). The article is structured as follows: the main results are stated in Section \ref{section:main}, the proofs of the results on state-constrained controllability and minimal time are respectively developed in Sections \ref{section:proofs1} and \ref{section:proofs2}, and we provide some perspectives for future research in Section \ref{section:perspectives}. 

\section{Main results \label{section:main}}

In the following, for $Y$ a vector in $\mathbb{R}^n$ and $\alpha \in \mathbb{R}$, we write $Y\geqslant \alpha$ and $Y>\alpha$ if all of the $n$ components of $Y$ are respectively greater or equal to $\alpha$ and greater than $\alpha$. Moreover, $| Y |$ is the usual Euclidean norm of $Y$ on $\mathbb{R}^n$ and $\max Y$ refers to the greatest component of $Y$. Finally, for $r\in \mathbb Z$ and $Q$ some open subset of $\mathbb R^{q}$ with $q\in\mathbb N^*$, $\mathrm{H}^r(Q)$ denotes the Sobolev space $\mathrm{W}^{r,2}(Q)$. 
In what follows, we will also consider the free evolution 

\begin{equation}
\label{eq:paralinfree}
\left \{
\begin{array}{l l l}
\partial_t \tilde Y - D \Delta \tilde Y &= A \tilde Y\mathbf &\text{in } \Omega_T, \\
\partial_{\overrightarrow{n}}\tilde Y&= 0 &\text{on } (0,T) \times \partial \Omega, \\
\tilde Y(0,\cdot)&=Y^0 (\cdot) &\text{in } \Omega.
\end{array}
\right .
\end{equation}
It is well-known that, for every $Y^0 \in \mathrm{L}^2(\Omega)^n$ and $U \in \mathrm{L}^2 (\omega_T)^m$, the Cauchy problem given by System \eqref{eq:paralin} admits a unique solution $Y \in \mathrm{L}^2((0,T),\mathrm{H}^1 (\Omega)^n) \cap \mathrm{H}^1((0,T),\mathrm{H}^{-1}(\Omega)^n)\hookrightarrow \mathrm{C}^0 ([0,T],\mathrm{L}^2(\Omega)^n)$.  If, in addition, $Y^0\in \mathrm{L}^\infty(\Omega)^n$ and $U\in \mathrm{L}^\infty(\omega_T)^m$, \changed{we also get a standard well-posedness $\mathrm{L}^{\infty}$ estimation:}

\begin{prpstn}[Well-posedness]
Let $T>0$. There exists $C(T) > 0$ such that, for any $Y^0 \in \mathrm{L}^{\infty}(\Omega)^n$ and any $U \in \mathrm{L}^{\infty}(\omega_T)^m$, the solution of \eqref{eq:paralin} with initial condition $Y^0$ and control $U$ satisfies
\begin{equation}
\| Y(t,\cdot) -  \tilde{Y}(t,\cdot) \|_{\mathrm{L}^{\infty} (\Omega_T)^n} \leqslant C(T) \| U \|_{\mathrm{L}^{\infty} (\omega_T)^m},
\label{eq:well}
\end{equation}
where $\tilde{Y}$ is the solution of the free System \eqref{eq:paralinfree} with same initial condition $Y^0$.
\label{prpstn:well}
\end{prpstn}

\changed{
\begin{proof}
Let $T$, $Y$ and $\tilde{Y}$ defined as in \ref{prpstn:well}, and let $Z = Y - \tilde{Y}$; $Z$ then satisfies the equation 
\begin{equation}
\label{eq:paralin-zero}
\left \{
\begin{array}{l l l}
\partial_t Z - D \Delta Z &= A Z + B U \mathbf{1}_{\omega} &\text{in } \Omega_T, \\
\partial_{\overrightarrow{n}} Z & = 0 &\text{on } (0,T) \times \partial \Omega, \\
Z(0,\cdot)&= 0 &\text{in } \Omega.
\end{array}
\right .
\end{equation}
System \eqref{eq:paralin-zero} is linear, so we have of course that the application $\mathcal{Z} : U \in \mathrm{L}^{\infty}(\omega_T)^m \mapsto  Z \in \mathrm{L}^{\infty}(\Omega_T)^n$ is linear. Moreover, by virtue of \cite[Chapter VII, Theorem 2.1]{ladyzhenskaia1988linear}, we have the following estimation on $Z$:
\begin{equation}
    \| Z \|_{\mathrm{L}^{\infty}(\Omega_T)^n} \leqslant C_1 + C_2 \| U \|_{\mathrm{L}^{\infty}(\omega_T)^m},
\end{equation}
where $C_1$ and $C_2$ depend only on $n$, $\Omega$, $T$ and $A$. Taking $U$ such that $\| U \|_{\mathrm{L}^{\infty}(\omega_T)^m} =1$, we obtain $\| Z \|_{\mathrm{L}^{\infty}(\Omega_T)^n} \leqslant C_1 + C_2$, so the linear map $\mathcal Z$ is also continuous, which yields estimation \eqref{eq:well} for all $Y^0 \in \mathrm{L}^{\infty}(\Omega)^n$ and $U \in \mathrm{L}^{\infty}(\omega_T)^m$, with $C(T) = C_1+C_2$.
\end{proof}
}


\subsection{Results on state-constrained controllability}

Before stating our main results, let us state a few preliminary properties satisfied by System \eqref{eq:paralin}. To preserve the nonnegativity of the trajectories of the uncontrolled System \eqref{eq:paralinfree}, we assume that the diffusion matrix $D$ is diagonal:
\begin{equation}
\forall i,j \in \{1,\dots, n\}, i \neq j \Rightarrow d_{ij}=0,
\label{eq:diagD}
\end{equation}
and that the coupling matrix $A$ is quasipositive, i.e.
\begin{equation}
\forall i,j \in \{ 1,\dots,n \}, i \neq j \Rightarrow a_{ij} \geqslant 0.
\label{eq:quasipos}
\end{equation}

Then, we have the following property:

\begin{prpstn}[{Positivity \cite[Lemma 1.1]{pierre_global_2010}}]
Assume \eqref{eq:diagD}. For any initial condition  $Y^0 \in \mathrm{L}^{\infty}(\Omega)^n$ such that $Y^0 \geqslant 0$, the corresponding solution $\tilde Y$ of \eqref{eq:paralinfree} satisfies $\tilde Y(t,\cdot)\geqslant 0$ for all time $t\geqslant 0$ if and only if \eqref{eq:quasipos} is satisfied. 
\label{prpstn:quasipos}
\end{prpstn}

Our goal is to control System \eqref{eq:paralin} with a positivity constraint on the state $Y$. Therefore, it seems reasonable to assume that free trajectories naturally stay nonnegative, hence assumptions \eqref{eq:diagD} and \eqref{eq:quasipos}.

We recall the classical notion of controllability to (free) trajectories:

\begin{dfntn}[Controllability to trajectories]
Let $T>0$. System \eqref{eq:paralin} is \textit{controllable to trajectories} in time $T$ if, for all $Y^0 \in \mathrm{L}^{2} (\Omega)^n$, and for all solution  $\tilde{Y}$ of the free system \eqref{eq:paralinfree} associated to another initial condition $\tilde{Y}^{0} \in \mathrm{L}^{2}(\Omega)^n$, there exists a control $U$ in $\mathrm{L}^{2} (\omega_T)^m$ such that the solution $Y$ of \eqref{eq:paralin} with initial condition $Y^0$ and control $U$ satisfies
\begin{equation}
Y(T,\cdot) = \tilde{Y} (T,\cdot).
\end{equation}
\end{dfntn}

The Laplace operator $-\Delta$ with Neumann boundary conditions admits a sequence of eigenvalues repeated with their multiplicity $(\lambda_p)_{p\in \mathbb{N}}$ such that 
\[
0 = \lambda_0 < \lambda_1 \leqslant \lambda_2 \leqslant \dots , \quad \lim_{p \rightarrow + \infty} \lambda_p + \infty.
\]
To each eigenvalue $\lambda_p$, we associate a corresponding normalized eigenvector $e_p\in \mathrm{H}^1(\Omega)$ in such a way that $\{e_p\}_{p\in\mathbb N}$ forms a Hilbert basis of $\mathrm{L}^2(\Omega)$. Notice that $e_p\in \mathrm{C}^\infty(\Omega)$ by elliptic regularity, since $\partial\Omega$ is smooth.
Given two matrices $\tilde A$ in $\mathcal{M}_{n} (\mathbb{R})$ and $\tilde B$ in $\mathcal{M}_{n,m} (\mathbb{R})$, we use the following notation for the \textit{Kalman matrix} associated to $\tilde{A}$ and $\tilde{B}$: 
\begin{equation}
\left [\tilde A | \tilde B \right ] = \left [ {\tilde A}^{n-1} \tilde{B} \; \big | \; {\tilde A}^{n-2} \tilde{B} \; \big | \; \dots \; \big | \; \tilde{B} \right ].
\end{equation}

Controllability to trajectories is ensured for the System \eqref{eq:paralin} under a Kalman-type condition:

\begin{prpstn}[{Controllability \cite{ammar-khodja_kalman_2008}}]
System \eqref{eq:paralin} is controllable to trajectories at any time $T$ if and only if, for all $p \in \mathbb{N}$, 
\begin{equation}
\mathrm{rank} \left [ (-\lambda_p D + A ) | B \right ] = n.
\label{eq:kalman}
\end{equation}
\label{prpstn:kalman}
\end{prpstn}

\begin{rmrk}
In the case $D=I_n$, condition \eqref{prpstn:kalman} simply becomes
\begin{equation}
\mathrm{rank} \left [ A  | B \right ] = n.
\label{eq:kalman2}
\end{equation}
\end{rmrk}

Our first two main results establish state-constrained controllability of the System \eqref{eq:paralin} under the aforementioned quasipositivity and controllability assumptions. 

\begin{thrm}[Case $D=I_n$]
Assume that $D=I_n$ and that $A$ and $B$ satisfy \eqref{eq:quasipos} and \eqref{eq:kalman2}. Assume moreover that the eigenvalues of $A$ all have a nonnegative real part.

Let $Y^0$, $Y^{\mathrm{f},0}$ in $\mathrm{L}^{\infty}(\Omega)^n$ and $Y^{\mathrm{f}}$ the solution of \eqref{eq:paralinfree} with initial condition $Y^{\mathrm{f},0}$. Assume that 
\begin{equation}
Y^0 \geqslant 0, Y^{\mathrm{f},0} \geqslant 0,
\label{eq:ssmptn-id1}
\end{equation} 
and that none of the components of $Y^0$ and $Y^{\mathrm{f},0}$ is a.e. zero on $\Omega$.

Then, there exists $T>0$ and $U \in \mathrm{C}^{\infty}_0(\omega_T)^m$ such that the solution $Y$ of \eqref{eq:paralin} with initial condition $Y^0$ and control $U$ satisfies 
\begin{equation}
Y(T,\cdot)=Y^{\mathrm{f}}(T,\cdot),
\label{eq:id1}
\end{equation}
and, for all $t$ in $[0,T]$,
\begin{equation}
Y(t,\cdot) \geqslant 0.
\label{eq:id2}
\end{equation}
\label{thrm:main-id}
\end{thrm}

\changed{\begin{rmrk}
The result above naturally still holds if we take $D = \alpha I_n$ for some $\alpha >0$, as it can be seen from its proof.
\end{rmrk}
}

\begin{thrm}[General case]
Assume that $A,B$ and $D$ satisfy \eqref{eq:ellipticity}, \eqref{eq:diagD}, \eqref{eq:quasipos} and \eqref{eq:kalman}. \changed{Let $Y^0$, $Y^{\mathrm{f},0}$ in $\mathrm{L}^{\infty}(\Omega)^n$ and $Y^{\mathrm{f}}$ the solution of \eqref{eq:paralinfree} with initial condition $Y^{\mathrm{f},0}$.} Assume that 
\[
Y^0 \geqslant 0, Y^{\mathrm{f},0} \geqslant 0.
\]
Then, for all $\varepsilon >0$, there exists $T>0$ and $U \in \mathrm{C}^{\infty}_0(\omega_T)^m$ such that the solution $Y$ of \eqref{eq:paralin} with initial condition $Y^0$ and control $U$ satisfies 
\begin{equation}
Y(T,\cdot)=Y^{\mathrm{f}}(T,\cdot),
\label{eq:main1}
\end{equation}
and, for all $t$ in $[0,T]$,
\begin{equation}
Y(t,\cdot) \geqslant - \varepsilon.
\label{eq:main2}
\end{equation}
\label{thrm:main}
\end{thrm}

\begin{rmrk}
Theorem \ref{thrm:main} can be indifferently stated with Neumann,  Dirichlet or even Robin boundary conditions in \eqref{eq:paralin}.

\changed{By contrast, Neumann conditions are an important assumption for Theorem \ref{thrm:main-id}, to ensure that the free trajectories of \eqref{eq:paralin}, after a well-chosen change of variables, converge to strictly positive (constant) steady states. Hence, Robin boundary conditions in Theorem \ref{thrm:main-id} could also be considered, as long as they make the free trajectories converge to positive steady states (possibly non-constant). 

On the other hand, since they prevent the solutions to be strictly positive, Dirichlet conditions appear to make the problem of controllability with nonnegative state more difficult and perhaps less relevant; the notion of approximate nonnegative controllability used in Theorem \ref{thrm:main} might then be better fitted to deal with such boundary conditions.
}  
\end{rmrk}

\changed{
\begin{rmrk}
Our proof enables to consider smooth controls $U \in \mathrm{C}^{\infty}_0(\omega_T)^m$ in the statement of Theorems \ref{thrm:main-id} and \ref{thrm:main}. The results still hold if we only have $\mathrm{L}^{\infty}_0(\omega_T)^m$ controls, \textit{e.g.} by replacing our control cost estimation \eqref{eq:cost} with a weaker one.
\end{rmrk}
}

\begin{rmrk}
Theorem \ref{thrm:main} establishes a form of ``approximate nonnegative controllability''. Moreover, it is sharp in the sense that exact nonnegative controllability does not hold in general, as illustrated by the following counterexample. Consider System \eqref{eq:paralin} with $d=1$, $\Omega = (0,1)$, $\omega$ any open subset of $\Omega$, $n=2$, $m=1$, $D=I_2$, $A = \begin{pmatrix} 0 & 1 \\ 0 & -1 \end{pmatrix}$ and $B = \begin{pmatrix} 0 \\ 1 \end{pmatrix}$, which leads to the following system:
\begin{equation}
\label{eq:example}
\left \{
\begin{array}{l l l}
\partial_t y_1 - \partial_{xx} y_1 &= y_2 &\text{in } (0,T) \times (0,1),\\
\partial_t y_2 - \partial_{xx} y_2 &= -y_2 + \mathbf{1}_{\omega} u &\text{in } (0,T) \times (0,1), \\
\partial_x Y (0) = \partial_x Y (1)  & = 0, \\
Y(0,\cdot)&=Y^0 (\cdot) &\text{in } (0,1).
\end{array}
\right .
\end{equation}
It is easy to see that System \eqref{eq:example} satisfy \eqref{eq:quasipos} and \eqref{eq:kalman}, so it is controllable with approximately nonnegative state by virtue of Theorem \ref{thrm:main}. Let us now attempt to require that the state remains exactly nonnegative. First, note that if $u=0$, $z = y_1+y_2$ satisfies the heat equation $\partial_t z - \partial_{xx} z = 0$, so for all $t\geqslant 0$,
\begin{equation} 
\int_0^1 z(t,x) \mathrm{d} x = \int_0^1 z(0,x) \mathrm{d} x.
\label{eq:example-2}
\end{equation} 
Now, we consider constant initial conditions $Y^0 = \begin{pmatrix} 3 \\ 1 \end{pmatrix}$ and $Y^{\mathrm{f},0} = \begin{pmatrix} 1 \\ 1 \end{pmatrix}$. The free trajectory $(y^{\mathrm{f}}_1,y^{\mathrm{f}}_2)$ starting at $Y^{\mathrm{f},0}$ 
verifies thanks to \eqref{eq:example-2} that $\int_0^1\left (y^{\mathrm{f}}_1(t,x)+y^{\mathrm{f}}_2(t,x)\right)\mathrm{d} x = 2$ for all $t\geqslant 0$. Moreover, Since $A$ is quasipositive in the sense of \eqref{eq:quasipos}, we have  $\int_0^1 y^{\mathrm{f}}_i(t,x)\geqslant 0$ for all $t\geqslant 0$ and $i\in\{1,2\}$, so $\int_0^1 y^{\mathrm{f}}_1(t,x) \mathrm{d} x \leqslant 2$ for all $t \geqslant 0$. On the other hand, since $y_2$ is nonnegative, the integral of $y_1$ on $[0,1]$ is nondecreasing over time, so the controlled trajectory $(y_1,y_2)$ satisfies $\int_0^1 y_1(t,x) \mathrm{d} x \geqslant 3$ for all $t \geqslant 0$, whichever the control. Therefore, it is impossible to control System \eqref{eq:example} from $Y^0$ to $Y^{\mathrm{f}}$.

This simple example highlights the fact that an actual gap exists between the notions of controllability with nonnegative state and approximately nonnegative state for coupled systems. In particular, Theorem \ref{thrm:main-id}  deals with a favorable case for which exact nonnegative controllability holds: when $D=I_n$ and the eigenvalues of $A$ have nonnegative real part. 

Let us additionally describe another situation in which exact nonnegative controllability holds between two trajectories, even with $D \neq I_n$. Assume that $\tilde{Y}$ and $Y^{\mathrm{f}}$ are globally bounded (if they are not, one can perform a change of variable $Y \mapsto e^{\lambda t}Y$ with $\lambda>0$ sufficiently large and apply the following to the new system). Then, let
 \begin{equation}
 \zeta = \min \left( \inf_{(\mathbb{R}_+ \times\Omega)^n} \tilde{Y} (t,x), \inf_{(\mathbb{R}_+ \times\Omega)^n} Y^{\mathrm{f}}(t,x) \right ).
 \label{eq:zeta}
 \end{equation}
\label{rmrk:problem}
If $\zeta > 0$, then one can replace $-\varepsilon$ with $\zeta - \varepsilon$ in \eqref{eq:main2} and conclude that \eqref{eq:paralin} is controllable between $\tilde{Y}$ and $Y^{\mathrm{f}}$ with nonnegative state, for $\varepsilon$ small enough. 

\changed{Thus, if $\tilde{Y}$ and $Y^{\mathrm{f}}$ are globally bounded and bounded from below by a positive constant, we recover nonnegative controllability.} The proof of this -- somewhat anecdotal -- result steadily follows the proof of Theorem \ref{thrm:main}, with $\zeta$ being added in the relevant inequalities. 
\end{rmrk}

The proofs of Theorems \ref{thrm:main-id} and \ref{thrm:main}, presented in Section \ref{section:proofs1}, are based on a ``staircase'' strategy, that has proven its efficiency for the study of state-constrained or control-constrained controllability \cite{loheac_minimal_2018,pighin_controllability_2018,pighin2019controllability,ruiz2020control}. The idea is to make small steps towards the target, following a path of trajectories such that the controlled trajectory stays always close to a nonnegative free trajectory, and therefore almost nonnegative (see Figures \ref{fig:staircase-id} and \ref{fig:staircase}). In the aforementioned references the steps trajectories are restricted to be connected steady states. The proof of Theorem \ref{thrm:main-id} features a change of variables that decouples the equations and allows the similar use of constant steady states. On the other hand, in the proof of Theorem \ref{thrm:main}, we relax this steady state assumption and follow a path of non-constant free trajectories. 

\subsection{Minimal time}

In this section, we do not assume anymore that $A,B$ and $D$ satisfy \eqref{eq:diagD}, \eqref{eq:quasipos} and \eqref{eq:kalman}. Our main result is the following:

\begin{thrm} Assume that assumption \eqref{eq:ellipticity} holds and that \textcolor{blue}{$\Omega\setminus\overline\omega$ contains a nonempty open ball}. Let $M>0$. Let $Y^0 \in \changed{\mathrm{L}^2(\Omega)^n}$, and let $Y^{\mathrm{f}} \in \mathrm{L}^2(\Omega_T)^n$ be a trajectory of System \eqref{eq:paralinfree} such that $Y^{\mathrm{f}}(0,\cdot) |_{\Omega \backslash \omega} \neq Y^0 |_{\Omega \backslash \omega}$. We define the minimal controllability time
\[
\bar{T}(Y^0,Y^{\mathrm{f}}) = \inf \left \{ T>0 / \exists u \in \mathrm{L}^2(\omega_T)^m \text{ s.t. } Y(0)=Y^0, Y(T,\cdot)=Y^{\mathrm{f}}(T,\cdot) \text{ and } \forall (t,x) \in \Omega_T, Y(t,x) \geqslant -M \right \},
\]
(with the convention $\inf \emptyset = + \infty$.) Then,  $\bar{T}(Y^0,Y^{\mathrm{f}}) > 0$.
\label{thrm:minimal}
\end{thrm}

\changed{
\begin{rmrk}
The controlled trajectories $Y(t,x)$ considered in the definition of $\bar{T}(Y^0,Y^{\mathrm{f}})$ are assumed to satisfy $Y(t,x) \geqslant -M$. In particular, if none of the controlled trajectories that go from $Y^0$ to $Y^{\mathrm{f}}$ admit a lower bound (in the $\mathrm{L}^{\infty}$ sense), then $\bar{T}(Y^0,Y^{\mathrm{f}}) = +\infty$.
\end{rmrk}
}

\begin{rmrk}
Under the assumptions \eqref{eq:diagD}, \eqref{eq:quasipos} and \eqref{eq:kalman} made in the previous section, Theorems \ref{thrm:main-id} and \ref{thrm:main} ensure moreover that $\bar{T}(Y^0,Y^{\mathrm{f}})  < + \infty$ either if $M<0$, or if $M=0$ and $D=I_n$ (or in the particular cases described in Remark \ref{rmrk:problem}). As mentioned in Remark \ref{rmrk:problem}, there exists cases for which  $\bar{T}(Y^0,Y^{\mathrm{f}}) = + \infty$ if $M=0$.
\end{rmrk}

\changed{
\begin{rmrk}
The proof of Theorem \ref{thrm:minimal} is based on a restriction of the solution in a ball \textcolor{blue}{strongly} included inside the domain $\Omega$, so the boundary conditions have little influence on the result; hence Theorem \ref{thrm:minimal} can be straightforwardly carried over for Dirichlet or Robin boundary conditions.
\end{rmrk}
}



\begin{rmrk}
The assumption $Y^{\mathrm{f}}(0,\cdot) |_{\Omega \backslash \omega} \neq Y^0 |_{\Omega \backslash \omega}$ corresponds to the ``interesting'' case where the control needs to act on regions over which it is not supported in order to reach the target trajectory. Therefore, we left out the case where $Y^{\mathrm{f}}(0,\cdot) |_{\Omega \backslash \omega}$ and $Y^0 |_{\Omega \backslash \omega}$ differ only on $\omega$, which however does not seem entirely trivial (notably, the strategy proposed in \cite[Remark 16]{loheac_minimal_2017} does not work). The positivity of the minimal time in that latter case may depend on whether $Y^{\mathrm{f}}$ is a constant steady state or a space-varying trajectory. 
\end{rmrk}

Theorem \ref{thrm:minimal} shows that relaxing the constraint $Y\geqslant 0$ to allow the controlled trajectory to be negative still implies the existence a minimal controllability time. This is not surprising: it has been numerically observed that, when there is no state constraint and as the time $T$ allowed to control the equation goes to zero, the control and the state tend to become highly oscillating \cite{boyer2013penalised,loheac_minimal_2017} and reach therefore very high absolute values. Hence, it is intuitively understandable that setting an unilateral constraint on the state restricts this behaviour and implies that $\bar{T}(Y^0,Y^{\mathrm{f}}) > 0$.

Theorem \ref{thrm:minimal} extends to general linear parabolic systems like \eqref{eq:paralin} the result stated in \cite[Theorem 4]{loheac_minimal_2017}, that establishes the existence of a positive minimal controllability time for the scalar heat equation. The proof, presented in Section \ref{section:proofs2}, relies on similar arguments as in \cite{loheac_minimal_2017}. 

\section{Proofs of Theorems \ref{thrm:main-id} and \ref{thrm:main} \label{section:proofs1}}

Before proving our results on state-constrained controllability, let us state a useful estimation on the cost of the control.

\begin{prpstn}[{Control cost}]\label{p:addreg}
Assume that \eqref{eq:kalman} holds. Let $Y^0 \in \mathrm{L}^{2}(\Omega)^n$ and let  $Y^{\mathrm{f}}$ be a trajectory of System \eqref{eq:paralinfree} associated to the initial condition $Y^{\mathrm{f},0} \in \mathrm{L}^{2}(\Omega)^n$. There exists a control $U \in \mathrm{C}^{\infty}_0(\omega_T)^m$ such that the corresponding solution $Y$ of \eqref{eq:paralin} satisfies $Y(T,\cdot) = Y^{\mathrm{f}}(T,\cdot)$, and satisfying moreover: for any $s\in \mathbb N$, there exists \changed{$C_s > 0$} such that 
\begin{equation}
\| U \|_{\mathrm{H}^s(\omega_T)^m} \leqslant \changed{C_s \exp \left ( \frac{C_s}{T} \right )} \| Y^0 - Y^{\mathrm{f},0} \|_{\mathrm{L}^2 (\Omega)^n}.
\label{eq:cost}
\end{equation}
\label{prpstn:cost}
\end{prpstn}
Proposition \ref{p:addreg} is classical but is not an immediate consequence of the results in \cite{ammar-khodja_kalman_2008}. For the sake of completeness, we provide a short proof, based on the strategy given in \cite[Theorem 4]{lissy2019internal}.
\begin{proof}
We consider the adjoint equation 

\begin{equation}\left\{ \begin{array}{lll} \label{heatadj}
\partial_t Z&=D^*\Delta Z +A^*Z&\text{in} \quad \Omega_T,\\
\partial_{\overrightarrow{n}} Z &= 0 &\text{on } (0,T) \times \partial \Omega, \\
Z(0,\cdot)&=Z^0&\text{in }\Omega.
\end{array} \right.
\end{equation}
First of all, we decompose the initial condition $Z^0$ in the Hilbert basis defined before \eqref{eq:kalman2}:  $$Z^0(x)=\sum_{k=1} Z^0_ke_k(x),\mbox{ }Z^0_k\in\mathbb R^n.$$ We can then decompose the  solution $Z$ of \eqref{heatadj} as 
$$Z(t,x)=\sum_{k=1}^\infty Z_k(t)e_k(x),$$ where $Z_k$ is the unique solution of the  ordinary differential system
\begin{equation}\label{fourode}
\left\{ \begin{array}{ll}
Z_k'&=(-\lambda_kD^*+A^*)Z_k,\\
Z_k(0)&=Z^0_k.
\end{array} \right .
\end{equation}
Let us recall the spectral  inequality for eigenfunctions of the Dirichlet-Laplace operator as obtained in the seminal paper by Gilles Lebeau and Enrique Zuazua \cite{MR1620510} (see also \cite{MR1743865}): for any non-empty open subset $\tilde\omega$ of $\Omega$, there exists $C>0$ such that for any $J\in \mathbb N^*$ and any $(a_1,\ldots a_J)\in \mathbb R^J$, we have 
\begin{equation}\label{lrscalaire}\sum_{k\leqslant J}|a_k|^2=\int_\Omega \left(\sum_{k\leqslant J} a_ke_k(x)\right)^2dx\leqslant C e^{C\sqrt{\lambda_J}} \int_{\tilde \omega} \left(\sum_{k\leqslant J} a_ke_k(x)\right)^2 \mathrm{d}x .\end{equation}
Writing \eqref{lrscalaire} for each component of $B^*Z$ and summing on $n$, we obtain that there exists $C>0$ such that, for all $t \in (O,T)$,
\begin{equation}\label{sie}\sum_{k\leqslant  J}\left|B^*Z_k(t)\right|^2\leqslant Ce^{C\sqrt{\lambda_J}} \int_{\tilde \omega} \left|\sum_{k\leqslant  J}B^*Z_k(t)e_k(x)\right|^2\mathrm{d}x.\end{equation}
Assume that $\tilde \omega$ is strongly included in $\omega$ and let $\varphi\in \mathrm{C}^\infty_0(\omega)$ be such that $\varphi=1$ on $\tilde \omega$. We can deduce from \eqref{sie} the inequality
\begin{equation}\label{lrscalaire2}\sum_{k\leqslant  J}|B^*Z_k(t)|^2\leqslant C e^{C\sqrt{\lambda_J}} \int_{\omega} \varphi(x) \left|\sum_{k\leqslant  J}B^*Z_k(t)e_k(x)\right|^2 \mathrm{d}x.\end{equation}
Integrating \eqref{lrscalaire2} between $T/4$ and $3T/4$, we obtain 
\begin{equation}\label{LR1}\int_{T/4}^{3T/4}\sum_{k\leqslant  J}|B^*Z_k(t)|^2dt\leqslant Ce^{C\sqrt{\lambda_J}}\int_{T/4}^{3T/4}\int_{\omega}\varphi(x) \left|\sum_{k\leqslant  J}B^*Z_k(t)e_k(x)\right|^2\mathrm{d}x\mathrm{d}t.\end{equation}

Now, we consider the system of ODEs \eqref{fourode}. Let $k\in \{1,\ldots,J\}$. Assumption \eqref{eq:kalman} implies that System \eqref{fourode} is observable on the time interval $(T/4,3T/4)$ and  we have the existence of some constant $C(\lambda_k)>0$ such that 
\begin{equation}\label{LR2}\left|Z_k\left(\frac{3T}{4}\right)\right|^2\leqslant C(\lambda_k) \int_{T/4}^{3T/4} ||B^*Z_k(t)||^2\mathrm{d}t.\end{equation}
Moreover, it is proved in \cite[Appendix]{lissy2019internal} that there exists $p_1,p_2\in\mathbb N$ (depending on $n$ but independent of $k$) such that \eqref{LR2} holds with
\begin{equation}\label{clk}C(\lambda_k)\leqslant C\left(1+\frac{1}{T^{p_1}}\right)\lambda_k^{p_2}.\end{equation}
Since $-\lambda_kA^*+D^*$ is dissipative for $k$ large enough,  there exists $C>0$ independent on $k$ and $T$ such that 
\begin{equation}
\label{eq:clk-bis}
||Z_k\left(T\right)||^2\leqslant C e^{CT}||Z_k\left(\frac{3T}{4}\right)||^2.
\end{equation}
Hence, restricting to the case $T\leqslant 1$ and combining \eqref{eq:clk-bis} together with \eqref{LR2} and \eqref{lrscalaire2}, we deduce that for another constant $C>0$,
\begin{equation}\label{obsLR}||Z_k\left(T\right)||^2\leqslant C\left(1+\frac{1}{T^{p_1}}\right) e^{C\sqrt{\lambda_J}}\int_{T/4}^{3T/4}\int_{\omega}\varphi(x)\left|\sum_{k\leqslant  J}B^*Z_k(t)e_k(x)\right|^2\mathrm{d}x\mathrm{d}t.
\end{equation}
Let $\psi\in \mathrm{C}^\infty_0(0,T)$ such that $\psi=1$ on $(T/4,3T/4)$. We deduce from \eqref{obsLR} that 
\begin{equation}\label{obsLR2}||Z_k\left(T\right)||^2\leqslant C\left(1+\frac{1}{T^{p_1}}\right)e^{C\sqrt{\lambda_J}} \int_{0}^{T}\int_{\omega}\varphi(x)\psi(t)\left|\sum_{k\leqslant  J}B^*Z_k(t)e_k(x)\right|^2\mathrm{d}x\mathrm{d}t.
\end{equation}
Inequality \eqref{obsLR2} is a low-frequency observability inequality for the solutions of \eqref{heatadj}. It is well-known that it is equivalent to a partial controllability result for the solutions of \eqref{eq:paralin}. 
More precisely, we consider as an initial condition 
$$\widehat Y^0=Y^0 - Y^{\mathrm{f},0}.$$
Then, we deduce that there exists $U_J\in \mathrm{L}^2(\omega_T)$, such that the corresponding solution $\widehat Y$ of \eqref{eq:paralin} with initial condition $\widehat Y^0$ satisfies that $\langle Y(T,\cdot), e_j \rangle =0$ for any $j\in \mathbb N$ with $j\leqslant J$. Moreover, following \cite[Proof of Proposition 2]{MR1620510}, it is possible to prove that one can choose $U_j$  in the smooth class $\mathrm{C}^\infty_0(\omega_T)$, in such a way that for any $s\in\mathbb N$, we have
\begin{equation}
\| U_J \|_{\mathrm{H}^s(Q_T)} \leqslant C\left(1+\frac{1}{T^{p_1}}\right)\frac{C}{\sqrt{T}}\left(1+\lambda_J^s\right) e^{C\sqrt{\lambda_J}} \| Y^0 - Y^{\mathrm{f},0} \|_{\mathrm{L}^2 (\Omega)^n}.
\label{eq:costJ}
\end{equation}
Hence, by applying the Lebeau-Robbiano strategy as described in \textcolor{blue}{\cite{MR2679651}}, we can create a control  $U\in \mathrm{C}^\infty_0(\omega_T)$ (that is alternating phases of  $\mathrm{C}^\infty$ controls with compact support and phases of dissipation) satisfying the estimation \eqref{eq:cost} such that the corresponding solution of \eqref{eq:paralin} with initial condition $\tilde Y^0$ and control $U$ satisfies $\widehat Y(T,\cdot)=0$. Then, by linearity, $Y=\ Y^{\mathrm{f}}+\widehat Y$ is a solution of \eqref{eq:paralin} (associated to the control $U$) with initial condition $Y^0$ and satisfies $Y(T,\cdot) = Y^{\mathrm{f}}(T,\cdot)$.

\end{proof}

Combining Propositions \ref{prpstn:well} and \ref{prpstn:cost}, with $s$ large enough such that $\mathrm{H}^s(Q_T)\hookrightarrow \mathrm{L}^\infty(Q_T)$, which is possible by Sobolev embedding, yields the following result, that features an estimation of the $\mathrm{L}^{\infty}$ distance between the free trajectory and the controlled trajectory by the $\mathrm{L}^2$ distance between the initial states of the free and the target trajectories. 

\begin{prpstn}
Let $Y^0 \in \mathrm{L}^{\infty}(\Omega)^n$, $\tilde{Y}$ the corresponding solution of \eqref{eq:paralinfree}, and $Y^{\mathrm{f}}$ a trajectory of \eqref{eq:paralinfree} with an initial condition $Y^{\mathrm{f},0}\in \mathrm{L}^{\infty}(\Omega)^n$. Then, for all $T>0$, there exists a control $U\in \mathrm{C}^\infty_0(\omega_T)$ and a constant $C(T)$ such that the solution $Y$ of \eqref{eq:paralin} with initial condition $Y^0$ and control $U$ satisfies $Y(T,\cdot) = Y^{\mathrm{f}} (T,\cdot)$ and
\begin{equation}
\forall t \in [0,T], \| Y(t,\cdot) - \tilde{Y}(t,\cdot) \|_{\mathrm{L}^{\infty}(\Omega)^n} \leqslant C(T) \| Y^0 - Y^{\mathrm{f},0} \|_{\mathrm{L}^2( \Omega)^n}.
\label{eq:lmm1}
\end{equation}
\label{prpstn:1}
\end{prpstn}

Proposition \ref{prpstn:1} is a key ingredient in the proofs of Theorems \ref{thrm:main-id} and \ref{thrm:main}.

\begin{proof}[Proof of Theorem \ref{thrm:main-id}] For $t\geqslant 0$ and $Y \in \mathbb{R}^N$, define
\begin{equation}
Z = e^{-tA} Y.
\label{eq:Z}
\end{equation}
Notice that $Y$ is a solution of \eqref{eq:paralin} if and only if $Z$ is solution of the following nonautonomous system:
\begin{equation}
\label{eq:paralin-id}
\left \{
\begin{array}{l l l}
\partial_t Z- \Delta Z &= e^{-tA} B U \mathbf{1}_{\omega} &\text{in } \Omega_T, \\
\partial_{\overrightarrow{n}} Z &= 0 &\text{on } (0,T) \times \partial \Omega, \\
Z(0,\cdot)&=Y^0(\cdot) &\text{in } \Omega.
\end{array}
\right .
\end{equation}

It is obvious that System \eqref{eq:paralin-id} is controllable if and only if System \eqref{eq:paralin} is controllable. Moreover, we can state an estimation similar to Proposition \ref{prpstn:1} for the trajectories of \eqref{eq:paralin-id}, that takes into account its time dependency:

\begin{lmm}
Let $T_0\geqslant 0$. Let $Z^0 \in \mathrm{L}^{\infty}(\Omega)^n$, and $\tilde{Z}$ the solution of
\begin{equation}
\label{eq:paralin-id-t}
\left \{
\begin{array}{l l l}
\partial_t Z- \Delta Z &= e^{-tA} B U \mathbf{1}_{\omega} &\text{in } (T_0,T_0+T) \times \Omega, \\
\partial_{\overrightarrow{n}} Z &= 0 &\text{on } (T_0,T_0+T) \times \partial \Omega, \\
Z(T_0,\cdot)&=Z^0(\cdot) &\text{in } \Omega,
\end{array}
\right .
\end{equation}
with initial condition $Z^0$ and $U=0$. Let $Z^{\mathrm{f}}$ be a trajectory of \eqref{eq:paralin-id-t} with another initial condition $Z^{\mathrm{f},0}$ and $U=0$. Then, for all $T>0$, there exists a control $U$ in $C^\infty_0(\omega_T)$ and $C(T)$ independent of $T_0$ such that the solution $Z$ of the system \eqref{eq:paralin-id-t} satisfies $Z(T_0+T,\cdot) = Z^{\mathrm{f}} (T_0+T,\cdot)$ and
\begin{equation}
\forall t \in [T_0,T_0+T], \| Z(t,\cdot) - \tilde{Z}(t,\cdot) \|_{\mathrm{L}^{\infty}(\Omega)^n} \leqslant C(T) \| Z^0 - Z^{\mathrm{f},0} \|_{\mathrm{L}^{2}( \Omega)^n}.
\label{eq:lmm-id1}
\end{equation}
\label{lmm:id1}
\end{lmm}

\begin{proof}
Let $t \in [T_0,T_0+T]$. Since System \eqref{eq:paralin} is autonomous, combining Proposition \ref{prpstn:1} and the definition of $Z$ gives the existence of a constant $C_1(T)$ independent of $T_0$ such that $\| e^{tA} ( Z(t,\cdot) - \tilde{Z}(t,\cdot) ) \|_{\mathrm{L}^{\infty}(\Omega)^n} \leqslant C_1(T) \| Z^0 - Z^{\mathrm{f},0} \|_{\mathrm{L}^{2}( \Omega)^n}$, so
\begin{equation}
\| Z(t,\cdot) - \tilde{Z}(t,\cdot) \|_{\mathrm{L}^{\infty}(\Omega)^n} \leqslant C(T) \vertiii{e^{-tA}} . \| Z^0 - Z^{\mathrm{f},0} \|_{\mathrm{L}^{2}( \Omega)^n}.
\end{equation}
Moreover, since the eigenvalues of $A$ all have nonnegative real part, it means that there exists $K$ independent of $t$ such that $\vertiii{e^{-tA}} \leqslant K$ for all $t\geqslant 0$, which yields \eqref{eq:lmm-id1} with $C(T) = K C_1 (T).$
\end{proof}

Let us highlight the fact that the estimation given by this lemma necessarily requires the assumption made on the eigenvalues of $A$. If we relax this assumption, we lose the independence of $C(T)$ with respect to $T_0$, which is is a key point of the following proof, for we will control System \eqref{eq:paralin-id} on a number of consecutive time intervals that depends itself on the value of $C(T)$. 

Now, let $Z^{\mathrm{f}}$ be the solution of \eqref{eq:paralin-id} with initial condition $Z^{\mathrm{f},0} = Y^{\mathrm{f},0}$ and no control. We are going to show the existence of $T>0$ and a control $U$ such that the solution $Z$ of \eqref{eq:paralin-id} with initial condition $Z^0=Y^0$ and control $U$ satisfies 
\begin{equation}
Z(T,\cdot) = Z^{\mathrm{f}} (T,\cdot)
\label{eq:zcontr}
\end{equation}
and, for all $(t,x) \in \Omega_T$,
\begin{equation}
Z(t,x) \geqslant 0.
\label{eq:zpos}
\end{equation}
As already noted above, it is clear that such a control $U$ is such that the solution $Y$ of \eqref{eq:paralin} with initial condition $Y^0$ and control $U$ satisfies
\begin{equation}
Y(T,\cdot) = Y^{\mathrm{f}} (T,\cdot).
\label{eq:ycontr}
\end{equation}
Moreover, for all $(t,x) \in \Omega_T$,
\begin{equation}
Y(t,x) = e^{tA} Z(t,x) \geqslant 0,
\label{eq:ypos}
\end{equation}
because of \eqref{eq:zpos} and the fact that the exponential of a quasipositive matrix has only nonnegative entries (if $A$ is quasipositive, write $A=P+\alpha I_n$ with $\alpha \in \mathbb{R}$ such that $P$ has only nonnegative entries, then it is clear that $e^P$ is nonnegative and so is $e^A=e^{\alpha I_n} e^P = e^{\alpha} e^P$ \cite{berman_nonnegative_1994}). 

When $U=0$, \eqref{eq:paralin-id} becomes a system of $n$ decoupled parabolic equations. Then, using a spectral expansion, we immediately have that the solutions $\tilde{Z}^0$ and $Z^{\mathrm{f}}$ starting respectively at $Z^0$ and $Z^{\mathrm{f},0}$ converge in norms $\mathrm{L}^{\infty}$ and $\mathrm{L}^2$ (and any other $\mathrm{L}^p$-norm) to 
\[
\bar{Z} = \frac{1}{| \Omega |} \int_{\Omega} Z^0 \quad \mathrm{and} \quad \bar{Z}^{\mathrm{f}} = \frac{1}{| \Omega |} \int_{\Omega} Z^{\mathrm{f},0},
\]
with $|\Omega| = \int_{\Omega} 1 \mathrm{d}x$. Assumption \eqref{eq:ssmptn-id1} in Theorem \ref{thrm:main-id} ensures that all the components of $\bar{Z}$ and $\bar{Z}^{\mathrm{f}}$ are positive. Therefore, there exists $\zeta>0$ such that $\bar{Z} \geqslant \zeta$ and $\bar{Z}^{\mathrm{f}} \geqslant \zeta$.

\begin{figure}
\begin{center}
\includegraphics[width=0.7\textwidth]{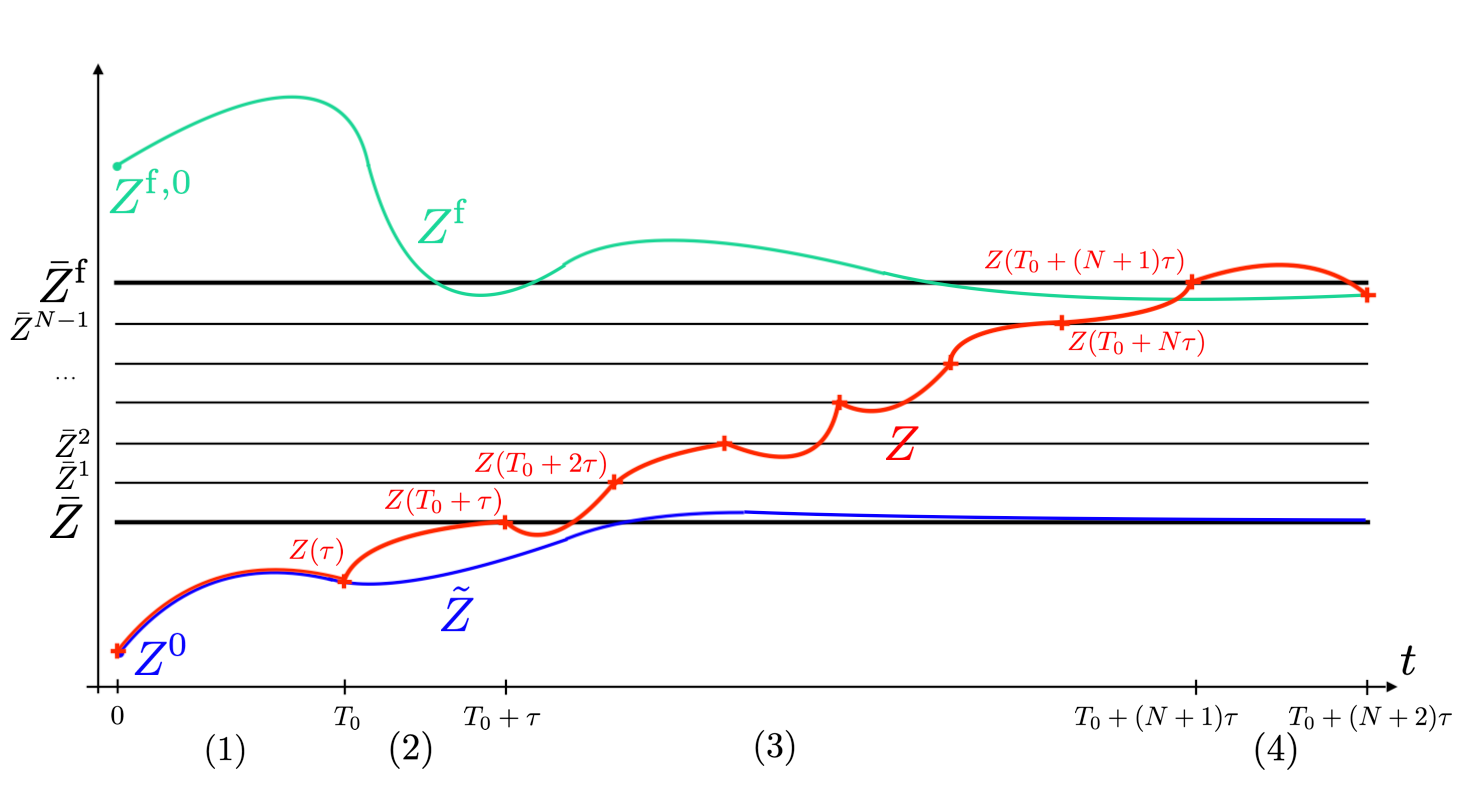}
\caption{Schematic representation of the staircase method used for the proof of Theorem \ref{thrm:main-id}.}
\label{fig:staircase-id}
\end{center}
\end{figure}

We are now in a position to build the trajectory $Z$ going from $Z^0$ to $Z^{\mathrm{f}}$. This will take several steps that are summarized on Figure \ref{fig:staircase-id}. Let $\delta >0$ and $\tau>0$.
\begin{enumerate}
\item We define a time $T_0$ (depending on $\delta$) such that, for all $t \geqslant T_0$,
\begin{equation}
\| \tilde{Z}_0 (t,\cdot) - \bar{Z} \|_{\mathrm{L}^{2}(\Omega)^n} \leqslant \delta \quad \mathrm{and} \quad \| Z^{\mathrm{f}} (t,\cdot) - \bar{Z}^{\mathrm{f}} \|_{\mathrm{L}^{2}(\Omega)^n} \leqslant \delta.
\label{eq:id-s1}
\end{equation}
The positivity of the initial conditions ensures that, with the control equal to $0$ on $[0,T_0]$,
\begin{equation}
\forall (t,x) \in [0,T_0] \times \Omega, Z(t,x) = \tilde{Z}_0 (t,x) \geqslant 0.
\end{equation}
\item  On the time interval $[T_0, T_0+\tau]$, let $V\in \mathrm{C}^{\infty}_0((T_0,T_0+\tau),\omega)^m$ be a control such that the solution ${Z}$ of 
\begin{equation*}
\left \{
\begin{array}{l l l}
\partial_t Z- \Delta Z  &= e^{-tA} B V \mathbf{1}_{\omega} &\text{in } \Omega_T, \\
\partial_{\overrightarrow{n}} Z &= 0 &\text{on } (T_0,T_0+\tau) \times \partial \Omega, \\
Z(T_0,\cdot)&=\tilde{Z}^0(T_0) &\text{in } \Omega.
\end{array}
\right .
\end{equation*}
satisfies $\tilde Z(T_0+\tau)=\bar{Z}$. Lemma \ref{lmm:id1} ensures that $V$ can be taken such that, for all $t \in [T_0, T_0+\tau]$, 
\begin{equation}
\| Z(t,\cdot) - \tilde{Z}_0(t,\cdot)  \|_{\mathrm{L}^{\infty}(\Omega)^n} \leqslant C(\tau) \| \tilde{Z}_0 (T_0) - \bar{Z} \|_{\mathrm{L}^{2}(\Omega)^n},
\end{equation}
and therefore, using \eqref{eq:id-s1}, for all $(t,x) \in [T_0, T_0+\tau] \times \Omega$,
\begin{equation}
Z(t,x) \geqslant \zeta - C(\tau) \delta.
\end{equation}
\item Let $M=\max (1, \max \bar{Z}_0, \max \bar{Z}^{\mathrm{f}}_0)$ and $N \in \mathbb{N}^*$ such that $\frac{M}{N} \leqslant \delta$. 

For $k \in \{ 0, \dots, N \}$, we define
\begin{equation}
\bar{Z}_k = (1-\frac{k}{N} )\bar{Z} + \frac{k}{N} \bar{Z}^{\mathrm{f},0},
\label{eq:k}
\end{equation}
(this way one has $\bar{Z}_0 = \bar{Z}$ and $\bar{Z}_N = \bar{Z}^{\mathrm{f}}$.) The $\bar{Z}^k$ define a sequence of constant steady states such that for all $k$ in $\{ 0, \dots, N-1 \}$, $\| \bar{Z}_{k+1} - \bar{Z}_{k} \|_{\mathrm{L}^{2}(\Omega)^n} \leqslant \delta$. 

Then, for each $k \in \{ 0, \dots, N-1 \},$ on the time interval $I_k=[T_0 + (k+1)\tau, T_0 + (k+2)\tau]$, we define a control $U_k$ in $\mathrm{C}^{\infty}_0(I_k \times \omega)^m$ such that the solution ${Z}$ of 
\begin{equation*}
\left \{
\begin{array}{l l l}
\partial_t Z- \Delta Z  &= e^{-tA} B U_k \mathbf{1}_{\omega} &\text{in } \Omega_T, \\
\partial_{\overrightarrow{n}} Z &= 0 &\text{on } I_k \times \partial \Omega, \\
Z(T_0+(k+1)\tau,\cdot)&=\bar{Z}_k &\text{in } \Omega.
\end{array}
\right .
\end{equation*}

satisfies $Z(T_0+(k+2)\tau)=\bar{Z}_{k+1}$. According to Lemma \ref{lmm:id1}, the control $U_k$ is such that one has, for all $t \in I_k$,
\begin{equation}
\| Z(t,\cdot) - \bar{Z}_{k} \|_{\mathrm{L}^{\infty}(\Omega)^n} \leqslant C(\tau) \| \bar{Z}_{k+1} - \bar{Z}_{k} \|_{\mathrm{L}^{2}(\Omega)^n},
\end{equation}
which means that
\begin{equation}
\forall (t,x) \in  I_k\times \Omega, Z(t,x) \geqslant \zeta - C(\tau) \delta.
\end{equation}
At the end of this step, we have reached the steady state $\bar{Z}_N=\bar{Z}^{\mathrm{f}}$.

\item On the time interval $[T_0+(N+1)\tau, T_0+(N+2)\tau]$, we define a control $W$ in $\mathrm{C}^{\infty}_0((T_0+(N+1)\tau, T_0+(N+2)\tau)\times\omega)^m$ such that the solution ${Z}$ of 
\begin{equation*}
\left \{
\begin{array}{l l l}
\partial_t Z- \Delta Z  &= e^{-tA} B W \mathbf{1}_{\omega} &\text{in } \Omega_T, \\
\partial_{\overrightarrow{n}} Z &= 0 &\text{on } (0,T) \times \partial \Omega, \\
Z(T_0+(N+1)\tau,\cdot)&=\bar{Z}^{\mathrm{f}} &\text{in } \Omega.
\end{array}
\right .
\end{equation*}
satisfies $Z(T_0+(N+2)\tau,\cdot)=Z_f(T_0+(N+2)\tau,\cdot)$. Lemma \ref{lmm:id1} ensures one more time that $W$ can be taken such that for all $t \in [T_0+(N+1)\tau, T_0+(N+2)\tau]$, 
\begin{equation}
\| Z(t,\cdot) - \bar{Z}^{\mathrm{f}} \|_{\mathrm{L}^{\infty}(\Omega)^n} \leqslant C(\tau) \| Z^{\mathrm{f}} (T_0+(N+1) \tau) - \bar{Z}^{\mathrm{f}} \|_{\mathrm{L}^{2}(\Omega)^n},
\end{equation}
which gives, using \eqref{eq:id-s1}: for all $(t,x) \in [T_0+(N+1) \tau,T_0+(N+2) \tau] \times \Omega$,
\begin{equation}
Z(t,x) \geqslant \zeta - C(\tau) \delta.
\end{equation}
\end{enumerate}
Overall, taking $T = T_0 + (N+2) \tau$,
\[
\delta = \frac{\zeta}{C(\tau)},
\]
which is possible since $\tau$ and $\zeta$ does not depend on $\delta$, and
\[
U(t) = \left \{ \begin{array}{l l}
0 & \text{on } [0,T_0), \\
V(t) & \text{on } [T_0, T_0 + \tau), \\
U_k(t) & \text{on } [T_0 + (k+1)\tau,T_0+ (k+2) \tau), k \in \{ 0, \dots N-1 \}, \\
W(t) & \text{on } [T_0+(N+1)\tau, T_0+(N+2)\tau],
\end{array} \right.
\]
we have found a control $U$ in $C^\infty_0(\omega_T)$ such that $Z$ satisfies \eqref{eq:zcontr} and \eqref{eq:zpos}.
\end{proof}

The proof of Theorem \ref{thrm:main} is again based on building a ``staircase'', made this time of non-constant trajectories.  Note that, to ensure that these trajectories do not go to far away from each other, we start the proof with a simple change of variable that makes the trajectories globally bounded. Like for the change of variable \eqref{eq:Z} in the proof of Theorem \ref{thrm:main-id}, this change of variable preserves the quasipositivity of the coupling matrix and the nonnegativity of the solutions.

\begin{proof}[Proof of Theorem \ref{thrm:main}]
Let $\tau>0$, $\varepsilon>0$,  $Y^0 \in \mathrm{L}^{\infty}(\Omega)^n$ such that $Y^0 \geqslant 0$, and $Y^{\mathrm{f}}$ in $\mathrm{L}^{\infty}(\mathbb{R}_+\times\Omega)^n$ a trajectory of \eqref{eq:paralin} associated to the initial condition $Y^{\mathrm{f},0} \geqslant 0$.

By means of a change of variable $Y \mapsto e^{\lambda t}Y$ with $\lambda>0$ sufficiently large (which is equivalent to changing $A$ into $A-\lambda I_n$, which does not affect the quasipositivity of the coupling matrix $A$), we can always assume that $A$ satisfies the following condition:
\begin{equation}
\forall \xi \in \mathbb{R}^N, \langle A \xi, \xi \rangle \leqslant 0.
\label{eq:bounded}
\end{equation}
We deduce from \eqref{eq:bounded} that the solution $Y$ of \eqref{eq:paralinfree} starting at any $Y^0$ is globally bounded. Indeed, take the scalar product of \eqref{eq:paralin} with $Y$ and integrate over $\Omega$:
\begin{equation*}
\int_{\Omega} \langle \partial_t Y, Y \rangle - \int_{\Omega} \langle D \Delta Y,Y \rangle = \int_{\Omega} \langle AY,Y \rangle,
\end{equation*}
which rewrites as 
\begin{equation*}
\frac{1}{2} \frac{d}{dt} \int_{\Omega} | Y |^2 = \int_{\Omega} \langle D Y,Y \rangle +  \int_{\Omega} \langle AY,Y \rangle \leqslant 0,
\end{equation*}
because of assumption \eqref{eq:bounded}, an integration by parts and assumption \eqref{eq:ellipticity}. Hence, for all $t>0$,
\begin{equation}\label{Ydis}
\int_{\Omega} | Y(t) |^2 \leqslant \int_{\Omega} | Y^0 |^2.
\end{equation}

\begin{figure}
\begin{center}
\includegraphics[width=0.7\textwidth]{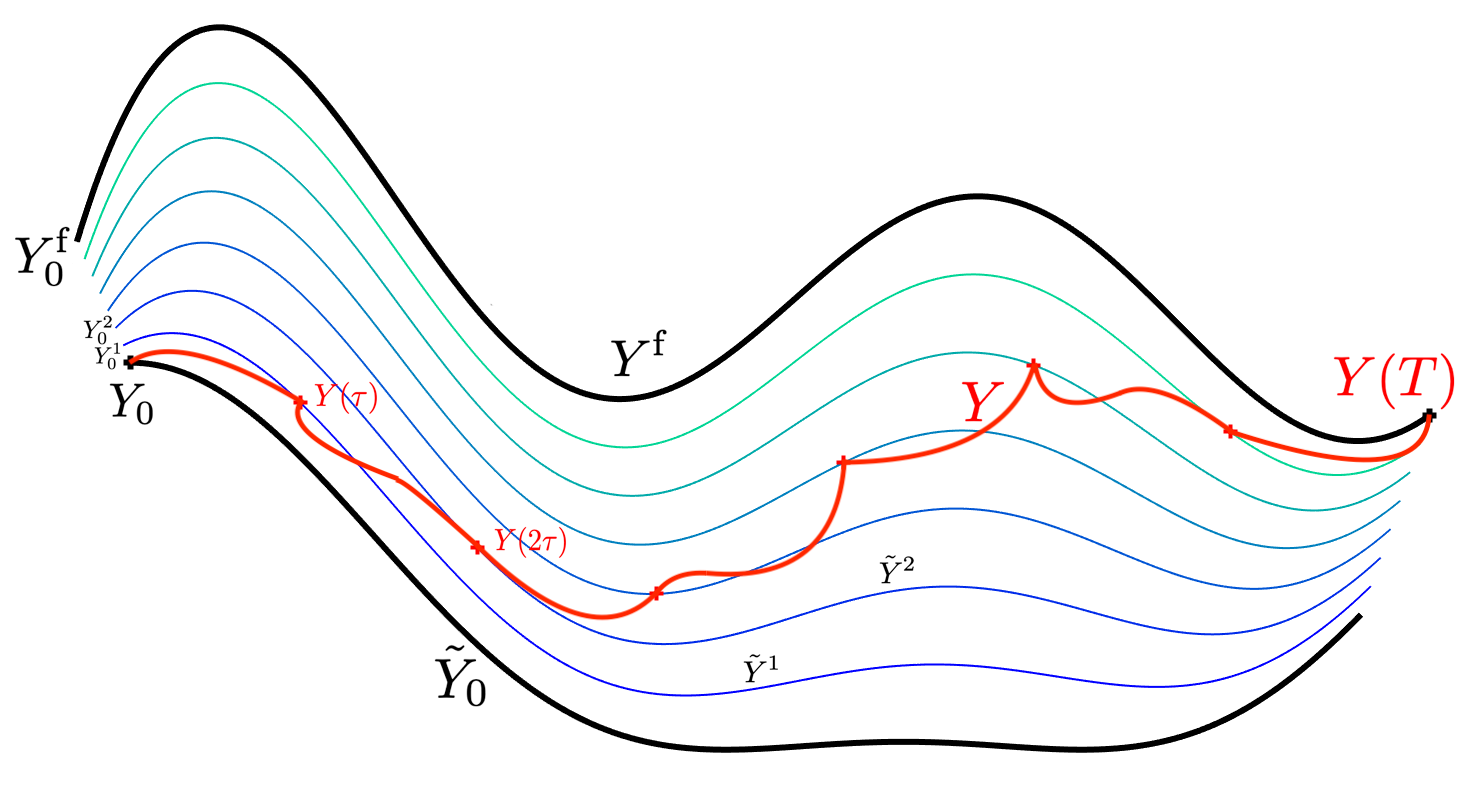}
\caption{Schematic representation of the staircase method used for the proof of Theorem \ref{thrm:main}.}
\label{fig:staircase}
\end{center}
\end{figure}

Let $\delta>0$. Let $N \in \mathbb{N}^*$ such that $\frac{1}{N} \| Y^0 - Y^{\mathrm{f},0} \|_{\mathrm{L}^{2}(\Omega)^n} \leqslant \delta$. For $k \in \{ 0, \dots, N \}$, we define 
\begin{equation}
Y^0_k = (1-\frac{k}{N} )Y^0 + \frac{k}{N} Y^{\mathrm{f},0},
\label{eq:yk}
\end{equation}
(this way one has $Y_0^0 = Y^0$ and $Y^0_N = Y^{\mathrm{f},0}$.) Let $\tilde{Y}_k$ be the solution of System \eqref{eq:paralinfree} with initial condition $Y^0_k$.
According to Proposition \ref{prpstn:quasipos}, $\tilde{Y}$ and $Y^{\mathrm{f}}$ are nonnegative. Let $M>0$ be such that $\| \tilde{Y} \|_{\mathrm{L}^{2}({\mathbb{R}_+ \times\Omega})} \leqslant M$ and $\| Y^{\mathrm{f}} \|_{\mathrm{L}^{2}({\mathbb{R}_+ \times\Omega})} \leqslant M$ ($M$ exists thanks to \eqref{Ydis}). Then, one has 
\[
\| \tilde{Y} -Y^{\mathrm{f}} \|_{\mathrm{L}^{2}({\mathbb{R}_+ \times\Omega})} \leqslant M,
\]
and it follows by linearity of System \eqref{eq:paralinfree} and the definition of the $\tilde{Y}_k$ that, for all $k \in \{ 0, \dots, N -1 \}$, one has
\begin{equation}
\| \tilde{Y}_{k+1} - \tilde{Y}_{k} \|_{\mathrm{L}^{2}({\mathbb{R}_+ \times\Omega})} \leqslant M \delta.
\label{eq:lmm2}
\end{equation}

Notice that, for all $k \in \{ 0, \dots, N \}$, $t\geqslant 0$ and $x \in \Omega$,
\[
\tilde{Y}_k(t,x) \geqslant 0.
\]

We now build the controlled trajectory $Y$ using the staircase strategy. The steps of the construction of $Y$ are represented on Figure \ref{fig:staircase}.

Let us start by steering the system from $Y^0$ to the trajectory $\tilde{Y}_1$. According to Proposition \ref{prpstn:1}, there exists a control $U_1 \in \mathrm{C}^\infty_0((0,\tau)\times\omega)^m$ such that the solution of \eqref{eq:paralin} with initial condition $Y^0$ and control $U_1$ satisfies 
\begin{equation}
Y(\tau,\cdot) = \tilde{Y}_1(\tau,\cdot),
\end{equation}
and, for all  $t \in [0,\tau]$,
\begin{equation}
\| Y(t,\cdot) - \tilde{Y} (t,\cdot) \|_{\mathrm{L}^{\infty}(\Omega)^n} \leqslant C(\tau) \| Y^0 - Y_1^0 \|_{\mathrm{L}^{2}(\Omega)^n}.
\end{equation}
Using \eqref{eq:lmm2}, we get
\begin{equation}
\forall t \in [0,\tau], Y(t,\cdot) \geqslant - C(\tau) M \delta.
\end{equation}

Next, let $k \in \{ 1, \dots, N-1 \}.$ We repeat the step above to steer the trajectory $Y$ from $\tilde{Y}_k(\cdot, k \tau)$ to $\tilde{Y}_{k+1} (\cdot, (k+1) \tau)$ in time $(k+1)\tau$ with a control $U_{k+1} \in \mathrm{C}^\infty_0((k\tau,(k+1)\tau)\times\omega)^m$ such that, for all $t \in [k \tau, (k+1) \tau]$, 
\begin{equation}
\| Y(t,\cdot) - \tilde{Y}_k (t,\cdot) \|_{\mathrm{L}^{\infty}(\Omega)^n} \leqslant C(\tau) \| \tilde{Y}_{k} (\cdot,k \tau) - \tilde{Y}_{k+1} (\cdot, k \tau) \|_{\mathrm{L}^{2}(\Omega)^n}.
\end{equation}
Using \eqref{eq:lmm2} again, we have 
\begin{equation}
\forall t \in [k \tau, (k+1) \tau], Y(t,\cdot) \geqslant - C(\tau) M \delta.
\end{equation}

Overall, let us set
\begin{equation}
\delta = \frac{\varepsilon}{M C(\tau)},
\label{eq:delta}
\end{equation}
$T=N\tau$, and $U$ the control defined on $[0,T]$ by $U(t,\cdot)=U_{k} (t,\cdot)$ if $t \in (k \tau, (k+1) \tau)$. The solution $Y$ of \eqref{eq:paralin} starting at $Y^0$ and with control $U$ in $C^\infty_0(\omega_T)$ satisfies \eqref{eq:main1} and \eqref{eq:main2}, which concludes the proof.
\end{proof}

\section{Proof of Theorem \ref{thrm:minimal} \label{section:proofs2}}

Our proof of the existence of a positive minimal time for controllability of System \eqref{eq:paralin} relies on proving the existence of such a minimal time for a scalar heat equation with a potential, a source term and boundary control. The arguments are inspired by those presented in the proof of \cite[Theorem 4.1]{loheac_minimal_2017}, which proves the same result for the standard heat equation.  

\begin{proof}
Let $M>0$ and $Y^0 \in \mathrm{L}^2(\Omega)^N$, and let $Y^{\mathrm{f}} \in \mathrm{L}^2(\Omega_T)^N$ be a trajectory of System \eqref{eq:paralin}. Assume that there exists $T>0$ and $U \in \mathrm{L}^2(\omega_T)^m$ such that the solution of \eqref{eq:paralin} starting at $Y^0$ with control $U$ satisfies $Y(T,\cdot) = Y^{\mathrm{f}}(T,\cdot)$ and 
\begin{equation}
\forall (t,x) \in \Omega_T, Y(t,x) \geqslant -M.
\label{eq:minimal-1}
\end{equation}

Assume without loss of generality that there exists an open ball $B$ contained in $\Omega \backslash \omega$ and such that the first components of $Y^0$ and $Y^{\mathrm{f}}(T,\cdot)$ differ on $B$, and consider the restriction to $B$ of the first equation of System \eqref{eq:paralin}:
\begin{equation}
\left \{
\begin{array}{l l}
\partial_t{y_1} - \Delta y = a_{11} y_1 + f(t,x) & \text{in } (0,T) \times B, \\
y_1(t,x) = v(t,x) & \text{on } (0,T) \times \partial B.
\end{array}
\right .
\label{eq:minimal-system}
\end{equation}
In \eqref{eq:minimal-system}, $f(t,x) = \sum_{j=2}^{N}  (a_{1j} + d_{1j} \Delta )y_j$ contains the coupling terms from System \eqref{eq:paralin}, and $v(t,x)$ is the trace of the solution $Y$ of \eqref{eq:paralin} with control $U$. \changed{Due to interior parabolic regularity results inside the domain $\Omega$, the solution $Y$ restricted to $B$  belongs to $\mathrm{L}^2((0,T),\mathrm{H}^2(B)^n)$, which notably ensures that $f \in \mathrm{L}^2((0,T)\times B)$.}  
Thus, \eqref{eq:minimal-system} can be seen as a scalar heat equation with a linear potential, a source term, and Neumann control on the whole boundary. Moreover, Assumption \eqref{eq:minimal-1} requires that the control $v$ in \eqref{eq:minimal-system} satisfies $v \geqslant -M$ at all times. 


Since we will only be considering System \eqref{eq:minimal-system} from now on, let us rename $y_1$ by $y$ and $a_{11}$ by $a$ to lighten notations. 

Let $y^0$ and $y^{\mathrm{f}}$ be the restrictions to $(0,T)\times B$ of the first component of $Y^0$ and $Y^{\mathrm{f}}$. Let also $y^{\mathrm{f},0} = y^{\mathrm{f}}(0,\cdot)$. We define as above the minimal controllability time for \eqref{eq:minimal-system} as
\begin{equation}
\begin{array}{l}
\bar{T}(y^0,y^{\mathrm{f}}) = \inf \big \{ T>0 / \exists v \in \mathrm{L}^2((0,T)\times \partial B) \\
\qquad \qquad \qquad \qquad \text{ s.t. } y(0,\cdot)=y^0, y(T,\cdot)=y^{\mathrm{f}}(T,\cdot) \text{ and } \forall (t,x) \in (0,T)\times \partial B, v(t,x) \geqslant -M \big \}. 
\end{array}
\end{equation}
Notice that by construction of \eqref{eq:minimal-system}, we have $\bar{T}(y^0,y^{\mathrm{f}}) \leqslant \bar{T}(Y^0,Y^{\mathrm{f}})$. Therefore, if \eqref{eq:minimal-system} has a positive minimal controllability time, then so does \eqref{eq:paralin}. 

Assume that $\bar{T}(y^0,y^{\mathrm{f}})=0$. Following the ideas used in \cite{loheac_minimal_2017} for proving the existence of a minimal time for the heat equation, we will study a spectral decomposition of the solution $y$ of \eqref{eq:minimal-system}. Consider the sequence of eigenvalues $(\lambda_n)_{n \in \mathbb{N}^*}$ and the associated sequence of eigenvectors $(p_n)$ of the following Sturm-Liouville problem on $[0,1]$:
\begin{equation}
\left \{ 
\begin{array}{l}
p_n''(r) + \frac{d-1}{r} p_n'(r) + a p_n(r) = - \lambda_n p_n (r), \\
p_n(1) = p_n'(0) = 0 \text{ and } p_n(0) \text{ s.t. } \omega_{d-1} \int_0^1 p_n^2 r^{d-1} \mathrm{d} r = 1.
\end{array}
\right.
\label{eq:sturm}
\end{equation}
In \eqref{eq:sturm}, $\omega_{d-1} = \int_{\partial B} \mathrm{d} \Gamma_x$. It is well-known that the sequence $(\lambda_n)$ is increasing and that $\lim_{n \rightarrow + \infty} \lambda_n = + \infty$.

We define $\varphi_n (x) = p_n(\| x \|)$ for $x \in B$, and $\alpha_n = p_n'(1)$, so that $\varphi_n$ satisfies the adjoint problem 
\begin{equation}
\left \{
\begin{array}{l l}
\Delta \varphi_n + a \varphi_n = - \lambda_n \varphi_n & \text{in } B, \\
\varphi_n (x) = 0, \nabla \varphi_n \cdot n(x) = \alpha_n & \text{on } \partial B.
\end{array}
\right.
\end{equation}
Notice that, thanks to the requirement made on $p_n(0)$ in \eqref{eq:sturm}, we have $\| \varphi_n \|_{\mathrm{L}^2(B)} = 1$ for all $n$ in $\mathbb{N}^*$. Moreover, straightforward computations (see \cite[Equation (18)]{loheac_minimal_2017}) give the identity
\begin{equation}
\| \varphi_n \|_{\mathrm{L}^2(B)} = \omega_{d-1} \frac{\alpha_n^2}{2 (\lambda_n+a)} (=1),
\label{eq:minimal-phi}
\end{equation}
that will be useful in the following.

Let $T>0$ and a control $v^T \in \mathrm{L}^2 ((0,T) \times \partial B)$ such that the solution $y$ of \eqref{eq:minimal-system} starting at $y^0$ with control $v^T$ reaches $y^{\mathrm{f}}$ in time $T$. For $n \in \mathbb{N}^*$ and $t \in (0,T)$, define $y_n = \int_B y(t,x) \varphi_n(x) \mathrm{d} x$. Then, we compute
\begin{align}
\dot{y}_n (t) & = \int_B (\Delta y(t,x) + a y(t,x) + f(t,x)) \varphi_n (x) \mathrm{d}x, \\
	& = \int_B y(t,x) \left ( \Delta \varphi_n + a \varphi_n \right ) \mathrm{d}x  -\int_{\partial B} y(t,x) \nabla \varphi_n \cdot n(x) \mathrm{d} \Gamma_x + \int_B f(t,x) \varphi_n (x) \mathrm{d} x, \\
	& = - \lambda_n y_n (t) - \alpha_n \int_{\partial B} v^T (t,x) \mathrm{d} \Gamma_x + \int_B f(t,x) \varphi_n (x) \mathrm{d} x,
\end{align}
which gives, after integrating on the time interval $(0,T)$,
\begin{equation}
y_n(T) = e^{-\lambda_n T} y_n(0) - \alpha_n \int_0^T e^{-\lambda_n (T-t)} \int_{\partial B} v^T \mathrm{d} \Gamma_x \mathrm{d} t + \int_0^T e^{-\lambda_n (T-t)} \int_B f \varphi_n \mathrm{d} x \mathrm{d} t,
\end{equation}
that we rewrite, using that $v^T = v^T + M - M$:
\begin{equation}
\begin{array}{l}
\displaystyle \int_0^T e^{-\lambda_n (T-t)} \int_{\partial B} (v^T+M) \mathrm{d} \Gamma_x \mathrm{d} t \\
\displaystyle \hspace{2em} = -\frac{1}{\alpha_n} \left ( y^{\mathrm{f}}_n (T) - e^{-\lambda_n T} y_n(0) - \int_0^T e^{-\lambda_n (T-t)} \int_B f \varphi_n \mathrm{d} x \mathrm{d} t \right ) + M \int_0^T e^{-\lambda_n (T-t)} \int_{\partial B} \mathrm{d} \Gamma_x \mathrm{d} t.
\end{array}
\label{eq:minimal-yn}
\end{equation}
Now, since $v^T + M \geqslant 0$, we have upper and lower bounds on the right-hand side of \eqref{eq:minimal-yn}, depending on the sign of $\lambda_n$:
\begin{equation}
\begin{array}{l}
\displaystyle \int_0^T \int_{\partial B} (v^T+M) \mathrm{d} \Gamma_x \mathrm{d} t  \\
\displaystyle \hspace{2em} \leqslant -\frac{1}{\alpha_n} \left ( y^{\mathrm{f}}_n (T) - e^{-\lambda_n T} y_n(0) - \int_0^T e^{-\lambda_n (T-t)} \int_B f \varphi_n \mathrm{d} x \mathrm{d} t \right ) + M \int_0^T e^{-\lambda_n (T-t)} \int_{\partial B} \mathrm{d} \Gamma_x \mathrm{d} t \\
\displaystyle \hspace{4em} \leqslant e^{-\lambda_n T} \int_0^T \int_{\partial B} (v^T+M) \mathrm{d} \Gamma_x \mathrm{d} t ,
\end{array}
\label{eq:minimal-yn-2-bis}
\end{equation}
if $\lambda_n \leqslant 0$, and
\begin{equation}
\begin{array}{l}
\displaystyle e^{-\lambda_n T} \int_0^T \int_{\partial B} (v^T+M) \mathrm{d} \Gamma_x \mathrm{d} t  \\
\displaystyle \hspace{2em} \leqslant -\frac{1}{\alpha_n} \left ( y^{\mathrm{f}}_n (T) - e^{-\lambda_n T} y_n(0) - \int_0^T e^{-\lambda_n (T-t)} \int_B f \varphi_n \mathrm{d} x \mathrm{d} t \right ) + M \int_0^T e^{-\lambda_n (T-t)} \int_{\partial B} \mathrm{d} \Gamma_x \mathrm{d} t \\
\displaystyle \hspace{4em} \leqslant \int_0^T \int_{\partial B} (v^T+M) \mathrm{d} \Gamma_x \mathrm{d} t,
\end{array}
\label{eq:minimal-yn-2}
\end{equation}
if $\lambda_n > 0$. Now, we want to take the limit when $T$ goes to 0 in \eqref{eq:minimal-yn-2-bis} and \eqref{eq:minimal-yn-2}. It is obvious that 
\begin{equation}
\lim_{T \rightarrow 0} M \int_0^T e^{-\lambda_n (T-t)} \int_{\partial B} \mathrm{d} \Gamma_x \mathrm{d} t = 0 \quad \text{ and } \quad \lim_{T \rightarrow 0} M \int_0^T \int_{\partial B} \mathrm{d} \Gamma_x \mathrm{d} t = 0.
\label{eq:minimal-limit1}
\end{equation}
Moreover, using the Cauchy-Schwarz inequality and the fact that $\| \phi_n \|_{\mathrm{L}^2(B)} =1$, we have the following bound:
\begin{equation} 
\left | \int_0^T e^{-\lambda_n (T-t)} \int_B f \varphi_n \mathrm{d} x \mathrm{d} t \right | \leqslant \| f \|_{\mathrm{L}^2((0,T) \times B)} \int_0^T e^{-\lambda_n (T-t)} \mathrm{d}t,
\end{equation}
and therefore 
\begin{equation}
\lim_{T \rightarrow 0} \int_0^T e^{-\lambda_n (T-t)} \int_B f \varphi_n \mathrm{d} x \mathrm{d} t  = 0.
\label{eq:minimal-limit2}
\end{equation}
Using \eqref{eq:minimal-limit1} and \eqref{eq:minimal-limit2} into \eqref{eq:minimal-yn-2-bis} and \eqref{eq:minimal-yn-2} when taking the limit yields
\begin{equation}
\lim_{T \rightarrow 0} \int_0^T \int_{\partial B} v^T \mathrm{d} \Gamma_x \mathrm{d} t = \frac{y_n^{\mathrm{f},0}  - y_n^0}{-\alpha_n}.
\label{eq:minimal-limit-end}
\end{equation}
Since the left-hand side of \eqref{eq:minimal-limit-end} does not depend on $n$, it means that there exists $\gamma \in \mathbb{R}$ such that for all $n \in \mathbb{N}^*$, 
\begin{equation}
y_n^0 = y_n^{\mathrm{f},0} + \alpha_n \gamma.
\label{eq:minimal-y}
\end{equation}

The next step is to show that $\gamma = 0$. Since $y^0 \in \mathrm{L}^{2} (B)$, we know that the series $\sum_{n=1}^{+\infty} | y_n^0 |^2$ converges. Plugging \eqref{eq:minimal-y} into this series yields
\begin{equation}
\sum_{n=1}^{+\infty} | y_n^0 |^2 = \sum_{n=1}^{+\infty} | y_n^{\mathrm{f},0} |^2 + \sum_{n=1}^{+\infty} \gamma \alpha_n ( 2 y_n^{\mathrm{f},0} + \alpha_n \gamma ).
\label{eq:minimal-series}
\end{equation}
The first sum on the right-hand side converges because $y^{\mathrm{f},0} \in \mathrm{L}^{2} (B)$. Therefore, the last sum is also finite, so 
\begin{equation}
\lim_{n \rightarrow \infty} \gamma \alpha_n ( 2 y_n^{\mathrm{f},0} + \alpha_n \gamma ) = 0.
\label{eq:minimal-end1}
\end{equation}
Then, notice that
\begin{equation*}
y_n^{\mathrm{f},0} = \int_{B} y^{\mathrm{f},0} \phi_n \mathrm{d} x = \frac{1}{\lambda_n + a} \int_{B} y^{\mathrm{f},0} (\lambda_n + a) \varphi_n \mathrm{d} x = \frac{1}{\lambda_n + a} \int_{B} y^{\mathrm{f},0} \Delta \varphi_n \mathrm{d} x = -\frac{\alpha_n}{\lambda_n + a} \int_{\partial B} y^1 \mathrm{d} \Gamma_x,
\end{equation*}
and use the identity \eqref{eq:minimal-phi} to obtain
\begin{equation}
\alpha_n y_n^{\mathrm{f},0} = \frac{2}{\omega_{d-1}} \int_{\partial B} y^{\mathrm{f},0} \mathrm{d} \Gamma_x,
\end{equation}
which does not depend on $n$. Since $\lim_{n \rightarrow +\infty} \alpha_n = + \infty$, we deduce that \eqref{eq:minimal-end1} holds only if $\gamma=0$.

This means that, for all $n \in \mathbb{N}$, $y_n^0 = y_n^{\mathrm{f},0}$. Therefore $y^0 = y^{\mathrm{f},0}$. By contraposition, this proves the theorem.
\end{proof}

\section{Discussion and open problems \label{section:perspectives}}

We have studied the problem of nonnegative controllability for coupled reaction-diffusion systems. Our results show that one can control such a system in large time to trajectories using the staircase method with approximately nonnegative state. Moreover, in the particular case where $D=I_n$ and $A$ only has eigenvalues with nonnegative real part, controllability in large time with nonnegative state holds. In a broader framework (less assumptions on $A$, $B$ and $D$), we also proved the existence of a positive minimal controllability time with whichever constraint of type $Y \geqslant -M$ with $M \geqslant 0$. 

We list a few remarks and open questions below.

\medskip

\textbf{Regularity of the control for minimal time.} In \cite{loheac_minimal_2017}, the authors show that the heat equation is controllable with nonnegative state constraint with a positive minimal time. Moreover, by considering a sequence of controls weakly converging in $\mathrm{L}^1$, they show that controllability in exactly the minimal time can be achieved with a Radon measure control. This result easily transposes to System \eqref{eq:paralin}. As stated in \cite{loheac_minimal_2017}, the question of whether the control in the minimal time can be more regular is still open.

\medskip

\textbf{Controllability with nonnegative state in the general case.} Remark \ref{rmrk:problem} displays an example of system showing that exact nonnegative controllability does not hold in general. Therefore, an interesting extension to Theorem \ref{thrm:main} would be to further discuss about the restrictions to be made on the initial condition and target state that could help recover exact nonnegative controllability. 

\medskip

\changed{\textbf{Non-autonomous systems.} Linear systems like \eqref{eq:paralin} with time-dependent matrices $A$, $B$ and $D$ (\textit{non-autonomous} systems) are also commonly considered and the question of their state-constrained controllability would be relevant. Controllability without a state constraint for such systems has been established in \cite{dupaix_generalization_2009} under a Silverman-Meadows-type condition. 

When adding a state constraint, our study suggests that estimations like \eqref{eq:ycontr} are crucial to establish controllability. As discussed after Lemma \ref{lmm:id1}, caution is required to guarantee that these estimations are uniform in time when the system has time-dependent coefficients. It is clearly not the case for all non-autonomous systems; hence finding conditions on $A$, $B$ and $D$ that allow controllability with non-negative state call for further investigation.}

\medskip

\changed{\textbf{Boundary control.}
Boundary control for coupled systems of parabolic equations is a difficult problem, and controllability even without a state constraint is not resolved as of today in the general case. The case $d=1$ and some particular cases when $d>1$ have been dealt with; we refer the reader to the survey paper \cite{ammar2011recent} and more recent advances made in \cite{benabdallah2014sharp,allonsius2020boundary}. A study of state-constrained controllability for these cases, potentially through straightforward adaptation of the staircase argument, would be an interesting continuation of this work. 
}

\medskip

\textbf{Nonlinear case.} A natural extension of this work would be the generalization of our results to semilinear parabolic systems. Let us do a short review of the state of the art for controllability and state-constrained controllability of such systems and give some perspectives on future research. Consider the following system

\begin{equation}
\label{eq:semilin}
\left \{
\begin{array}{l l l}
\partial_t y - \Delta y &= f(y) + B u \mathbf{1}_{\omega} &\text{in } \Omega_T, \\
\partial_{\overrightarrow{n}} y &= 0 &\text{on } (0,T) \times \partial \Omega, \\
y(0,\cdot)&= y^0 (\cdot) &\text{in } \Omega,
\end{array}
\right .
\end{equation}
with $f : \mathbb{R} \rightarrow \mathbb{R}$ a locally Lipschitz-continuous function satisfying $f(0) = 0$ and the following properties:
\begin{equation}
     \exists C > 0, \forall s \in \mathbb{R}, | f'(s) | \leqslant C (1+ |s|^{1+4/d}),
\end{equation}
and
\begin{equation}
    \exists \alpha > 0, \frac{f(s)}{| s | \log^{\alpha}\left ( 1+|s| \right )} \underset{|s| \rightarrow +\infty}{\longrightarrow} 0.
    \label{eq:alpha}
\end{equation}
Due to the nonlinearity $f$, in absence of a control $u$, the state can exhibit blowup in finite time. The controllability properties of System \eqref{eq:semilin} depend on the value of $\alpha$ in \eqref{eq:alpha}:
\begin{itemize}
    \item for $\alpha\leqslant 3/2$, \eqref{eq:semilin} is controllable to trajectories in arbitrary time \cite[Theorem 1.2]{fernandez-cara_null_2000},
    \item for $\alpha >2$, \eqref{eq:semilin} might fail to be controllable and blow up in finite time \cite[Theorem 1.1]{fernandez-cara_null_2000},
    \item for $3/2 < \alpha \leqslant 2$, assuming $s f(s) > 0$ for any $s \neq 0$, \eqref{eq:semilin} is null-controllable in large time \cite[Theorem 2.5]{balch_global_2018}.
\end{itemize}

As for nonnegative-state controllability, first results have been stated in \cite{nunez2019controllability} in two particular cases:
\begin{thrm}[{\cite[Theorems 1.1 and 1.2]{nunez2019controllability}}]
\begin{enumerate}
\item (Steady-state controllability). Let $y_0$ and $y_1$ in $\mathrm{L}^{\infty}(\Omega)$ be two positive steady states of \eqref{eq:semilin}. Assume $y_0$ and $y_1$ are connected, \textit{i.e.} there exists a continuous map $\gamma = [0,1] \rightarrow \mathrm{L}^{\infty}(\Omega)$ such that $\gamma(0) = y_0$ and $\gamma(1) =1$. Moreover, assume that for all $s \in [0,1]$, $\gamma(s) > 0$.

Then, there exists a time $T>0$ and a control $u \in \mathrm{L}^{\infty}((0,T)\times \Gamma)$ such that the solution $y$ of \eqref{eq:semilin} with initial condition $y_0$ satisfies $y(T) = y_1$ and, for all $t$ in $(0,T)$, $y(t) \geqslant 0$. 
\item (Controllability in the dissipative case). In the dissipative case ($sf(s)\geqslant 0$ for all $s \in \mathbb{R}$), System \eqref{eq:semilin} is controllable to trajectories in large time with nonnegative state. 
\end{enumerate}
\label{thrm:nunez}
\end{thrm}
In both cases, there exists a positive minimal time. Let us also mention \cite{pighin_controllability_2018}, in which are shown similar results for a semilinear equation with \textit{boundary} control, and also that state-constrained controllability fails outside of these two particular cases.

To extend the results of Theorem \ref{thrm:nunez} to more general nonlinearity $f$ or arbitrary initial and target data, the main challenge in the semilinear case compared to the linear case is that, even in a favourable case where the nonlinearity $f$ is globally Lipschitz, the staircase method does not work anymore, because the trajectories might move away from each other exponentially in time. Therefore, the small fixed-size steps of the staircase do not ensure that the controlled trajectory will eventually reach the target trajectory. We even conjecture that this type of behaviour might make state-constrained controllability fail, and are conducting research to find a counterexample. 

Finally, for a semilinear system of coupled equations,
\begin{equation}
\label{eq:semilin-syst}
\left \{
\begin{array}{l l l}
\partial_t Y - D \Delta Y &= f(Y) + B u \mathbf{1}_{\omega} &\text{in } \Omega_T, \\
\partial_{\overrightarrow{n}} Y &= 0 &\text{on } (0,T) \times \partial \Omega, \\
Y(0,\cdot)&= Y^0 (\cdot) &\text{in } \Omega,
\end{array}
\right .
\end{equation}
with $Y \in \mathbb{R}^N$, little is known on global controllability to trajectories. Moreover, for scalar equations, most state-constrained controllability results rely to some extent on the maximum principle, which does not hold for coupled systems like \eqref{eq:semilin-syst}. Hence the question of state-constrained controllability for these systems remains largely open. 

\section*{Acknowledgements}
The authors would like to thank Ariane Trescases for having pointed out a mistake in the proof of Theorem \ref{thrm:minimal}.
\bibliography{biblio2}
\bibliographystyle{IEEEtran}

\end{document}